\documentclass[11pt,a4paper]{article}
\usepackage[T1]{fontenc}
\usepackage[utf8]{inputenc}
\usepackage{authblk}

\usepackage{moreverb}
\usepackage[english]{babel}

\newcommand{\ds}{\displaystyle}

\newcommand{\E}{\mathbb{E}}
\newcommand{\V}{\mathbb{V}}

\usepackage{amsmath} 
\usepackage{amssymb} 
\usepackage{subcaption}
\usepackage{mathrsfs} 
\usepackage{fancyhdr} 
\usepackage{graphicx} 
\usepackage[T1]{fontenc} 
\usepackage[utf8]{inputenc} 
\usepackage{times} 
\usepackage{eurosym} 
\usepackage{color} 
\usepackage{multicol} 
\usepackage{ulem} 
\usepackage{makeidx}

\usepackage{caption}
\usepackage[]{natbib}
\usepackage[colorlinks,citecolor=blue]{hyperref}

\date{}

\begin{document}


\begin{center}
{\huge
	{A fully objective Bayesian approach for the Behrens-Fisher problem using historical studies}
}
\bigskip
\bigskip

\begin{tabular}{ccc}
\large{Antoine Barbieri}$^{1,2}$ & \large{Jean-Michel Marin}$^{1}$ & \large{Karine Florin}$^{3}$ \\ 
\footnotesize{(Antoine.Barbieri@umontpellier.fr)} & \footnotesize{(Jean-Michel.Marin@umontpellier.fr)} & \footnotesize{(Karine.Florin@sanofi.fr)} \\ 
\end{tabular} 
\bigskip

{\it
$^{1}$ {Institut Montpelliérain Alexander Grothendieck (IMAG), Université de Montpellier }\\ 
$^{2}$ {Institut régional du Cancer Montpellier - Val d’Aurelle (ICM), Unité de Biométrie }\\
$^{3}$ {Unité de Biostatistiques, Sanofi, Montpellier, France}
}
\end{center}

\bigskip

\begin{abstract}
For in vivo research experiments with small sample sizes and available historical data, we propose a sequential Bayesian method for the Behrens-Fisher problem. We consider it as a model choice question with two models in competition: one for which the two expectations are equal and one for which they are different. The choice between the two models is performed through a Bayesian analysis, based on a robust choice of combined objective and subjective priors, set 
on the parameters space and on the models space. Three steps are necessary to evaluate the posterior probability of each model using two historical datasets similar to the one of interest. Starting from the Jeffreys prior, a posterior using a first historical dataset is deduced and allows to calibrate the Normal-Gamma informative priors for the second historical dataset analysis, in addition to a uniform prior on the model space. From this second step, a new posterior on the parameter space and the models space can be used as the objective informative prior for the last Bayesian analysis. Bayesian and frequentist methods have been compared on simulated and real data. In accordance with FDA recommendations, control of type I and type II error rates has been 
evaluated. The proposed method controls them even if the historical experiments are not completely similar to the one of interest.
\end{abstract}


\section{Introduction}
 
In preclinical research, experiments are routinely performed using exactly the same protocol under the same conditions. Historical data are often available with the same control, the same reference product or similar products to the one tested. 
A specific feature of {\it in vivo} pharmacology experiments is the small sample size per experiment, due to ethical considerations \citep{landis_call_2012}.
Frequently, the objective of such experiments is to evaluate a treatment effect in comparison to a control (placebo or reference compound).  

Before dealing with the experimental context used to evaluate the effects of a compound at different doses or to compare different treatments, we focused on a simple study design of an experimental treatment group compared to a control group in parallel. 
The measurements are supposed to be normally distributed with expectation $\mu_t$ when the experimental treatment is administered and $\mu_c$ if not. 
The variances of the two Gaussian distributions are unknown and assumed to be different. The aim is to discriminate between the two hypotheses
\begin{equation}
\left\{
    \begin{array}{lll}
        \mathcal{H}_{0}: & \mu_{c}=\mu_{t}=\mu & \\
        \mathcal{H}_{1}: & \mu_{c}\neq\mu_{t} & ,
    \end{array}
\right. \label{eq_hyp}
\end{equation}
which is known as the Behrens-fisher problem \citep{kim_behrens-fisher_1998}.

In the case of hypothesis testing, the frequentist methods are routinely applied using the p-value interpretation. These methods do not take into account the historical data and are poorly adapted for small sample sizes.
The frequentist approach for solving this problem is the Student test with Satterthwaite correction.
However, the small sample size causes a lack of information in the statistical analysis, hence the idea of including some prior information in the inferential process, to increase the discriminative power of the procedure. 
With the frequentist method, a crude way to take into account historical datasets is to pool the experiments. One can test if the treatment effect is significant after pooling the data.
In contrast, the Bayesian paradigm is considered to be a well-suited method for handling small samples, in particular since it allows taking into account prior information.
In the Bayesian paradigm, two directions can be adopted:
\begin{itemize}
\item the first one is to use the historical datasets to determine informative prior distributions;
\item the second one is to use a hierarchical model which assumes a distribution across studies with a parameter that controls the variation.
\end{itemize}
There are numerous studies on the second case \citep{neuenschwander_summarizing_2010,hobbs_hierarchical_2011,viele_use_2014,schmidli_robust_2014}, which is an interesting and natural modeling of the problem. 
However, the results are very sensitive to the prior distribution on the between-experiments dispersion parameter \citep{gelman_prior_2006}. 
In this work, we investigated the first case with constraint on the sample size. 
In this Bayesian situation, work has previously been performed on sample size calculation \citep{whitehead_bayesian_2008,whitehead_bayesian_2015}, but that won't be addressed in our work.

For a given prior distribution, a typical procedure to discriminate between $\mathcal{H}_0$ and $\mathcal{H}_1$ is to construct a Bayesian credible interval of level $1-\alpha$ ($\alpha$ typically equal to $5\%$) on $\mu_c-\mu_t$ and to decide in favour of $\mathcal{H}_0$ if 0 belongs to the credible interval and $\mathcal{H}_1$ if not \citep{kim_behrens-fisher_1998,ghosh_behrens-fisher_2001}.
This is not the fully Bayesian approach.
In the Bayesian paradigm, a prior distribution on the space of hypotheses must be defined and then their posterior distribution calculated \citep{marin_bayesian_2013,kass_bayes_1995}. 
Indeed, each hypothesis is associated with a model: the model $\mathcal{M}_0$ under which $\mu_c=\mu_t$ and the model $\mathcal{M}_1$ under which $\mu_c\neq \mu_t$.
Then, prior probability of model $\mathcal{M}_0$, denoted $\Pr\left(\mathcal{M}_0\right)$, directly related to $\mathcal{M}_1$ (with $\Pr\left(\mathcal{M}_1\right)=1-\Pr\left(\mathcal{M}_0\right)$) is set and the posterior probability of $\mathcal{M}_0$ is calculated.
Using this last probability and a loss function, the experimenter can decide whether to reject or not $\mathcal{H}_0$.

There are very few works on a fully Bayesian approach for the Behrens-Fisher problem. To our knowledge, there is only one proposal of \citet{moreno_default_1999} in which intrinsic and fractional Bayes factors are calculated to avoid the improper prior difficulty. Here, we propose instead to use some historical data.
But, it is not easy to calibrate an objective prior distribution using historical datasets such that the model posterior probabilities are correctly defined.
In this work, we present a solution in three steps involving two similar historical datasets used sequentially to calibrate the prior.
In this work, we propose a solution in three steps involving two similar historical datasets used sequentially to calibrate the prior.
The first allows building an informative prior on the model parameters from an improper and non informative prior. 
The informative priors on the model and parameter spaces are calibrated at the completion of the second step using the proper informative prior on parameter from the first step and the second considered historical dataset.
The last step returns the posterior probability of each model, given the informative prior from the second step, allowing to make the statistical inference on the experiment of interest.
The plan of the paper is as follows.
In section 2, the methodology is introduced and the proposed Bayesian method is detailed. 
Section 3 then presents a simulation study to compare our method and the frequentist methods.
Finally, section 4 applies both approaches on real data from \textit{in vivo} behavioral pharmacology experiments.

\section{Methods}

\subsection{Recap on Bayesian model choice}

Let $y$ denote the observed dataset and $J$ the number of models in competition.
In the parametric Bayesian paradigm, model $j\in\{0,\ldots,J-1\}$, denoted $\mathcal{M}_{j}$, is defined with
\begin{itemize}
\item a likelihood function $\ell_j(\theta_j|y)$ with unknown parameter $\theta_j\in\Theta_j$;
\item a prior distribution denoted by $\pi_j(\theta_j)$ on the parameter space of $\mathcal{M}_{j}$.
\end{itemize}
The Bayesian model choice is done according to the model posterior probabilities \citep{robert_bayesian_2007,scott_exploration_2005}.
Given prior probabilities on the models $\Pr\left(\mathcal{M}_j\right)$, the posterior distribution on the model space is deduced by
$$
\Pr\left(\mathcal{M}_j|y\right)\propto \left\{\Pr\left(\mathcal{M}_j\right)\int_{\Theta_j}\ell_j(\theta_j|y)\pi_j(\theta_j)d\theta_j\right\}\,.
$$
With $\Pr(y|\mathcal{M}_j)=\int_{\Theta_j}\ell_j(\theta_j|y)\pi_j(\theta_j)d\theta_j$, the posterior probabilities become
\begin{equation} 
\label{eq_posterior}
\Pr\left(\mathcal{M}_{i}\mid y\right) = \frac{\Pr\left(\mathcal{M}_{i}\right) \Pr\left(y\mid \mathcal{M}_{i}\right)}{\sum_{j=0}^{J-1}\Pr\left(\mathcal{M}_{j}\right) \Pr\left(y\mid \mathcal{M}_{j}\right)}.
\end{equation}
The key quantities in equation \eqref{eq_posterior} are the $\Pr\left(y\mid \mathcal{M}_{j}\right)$ which are called marginal or integrated likelihoods.
In absence of a loss function, we select the model one with the largest posterior probability.

\subsection{The model}

A set of two independent and identically distributed (\textit{iid}) samples $c$ and $t$ are considered. The first one associated with the control group is assumed to be normally distributed of size $n_c$ with expectation $\mu_c$ and variance $\sigma_c^2$. 
The second one of size $n_t$ is associated with the treated group and assumed to be normally distributed with expectation $\mu_t$ and variance $\sigma_t^2$. 
We focused on the Behrens-Fisher problem where the aim is to test \eqref{eq_hyp}. 
In the fully Bayesian paradigm presented in the previous section, this hypothesis testing is equivalent to a model choice problem between the two following models:
\begin{itemize}
\item $\mathcal{M}_0$ under which $c\sim_{iid}\mathcal{N}(\mu,\sigma_{0,c}^2)$ and $t\sim_{iid}\mathcal{N}(\mu,\sigma_{0,t}^2)$,
\item $\mathcal{M}_1$ under which $c\sim_{iid}\mathcal{N}(\mu_c,\sigma_{1,c}^2)$ and $t\sim_{iid}\mathcal{N}(\mu_t,\sigma_{1,t}^2)$.
\end{itemize}
Let $\theta_0=(\mu,\sigma_{0,c}^2,\sigma_{0,t}^2)$, $\theta_1=(\mu_c,\mu_t,\sigma_{1,c}^2,\sigma_{1,t}^2)$
and $y=(c,t)$.  
According to model $j$, the posterior distribution of $\theta_j$ is:
$$
\pi_j\left(\theta_{j}\vert y\right)\propto\ell_j\left(\theta_{j}| y\right)\pi\left(\theta_{j}\right)\,.
$$
The aim of our work is to introduce a methodology that uses historical datasets to calibrate objective prior distributions on $\theta_0$ and $\theta_1$ and also on the model index.
We assume to have at disposal two previous datasets from similar experiments. 

\subsection{Our objective Bayesian approach}

Typically, to solve the Behrens-Fisher problem, the practitioners use the following Bayesian approach:
\begin{itemize}
\item to forget the model choice aspect of the problem;
\item to construct a credible interval on the difference between the two expectation;
\item to reject $H_0$ if 0 does not belong to the credible interval.
\end{itemize}
In this case, the Bayesian paradigm is just used to construct a credible interval. The prior distribution on the two hypotheses, i.e. the prior probabilities $\Pr(\mathcal{M}_0)$ and $\Pr(\mathcal{M}_1)$, do not exist. 
In the contrast, a fully Bayesian approach considers prior probabilities and deduces posterior probabilities on both hypotheses.
The benefit of using such a procedure is to weight each hypothesis according to some knowledge such as the data. 
Using these weights, the experimenter can then decide whether to reject $H_0$ or not. 
The setting of the prior distribution can be controversial, especially if it is based on expert opinion.
However, we show below how to calibrate the prior in presence of historical data.
Thus, a robust choice taking advantage of historical data to build an objective prior is proposed: starting from a non informative prior, the historical data are used to deduce an informative proper prior.

Priors must be defined for $\theta_0$, $\theta_1$ and also for the probabilities  $\Pr(\mathcal{M}_0)$ and $\Pr(\mathcal{M}_1)$.
Our proposal is to use a posterior distribution based on the historical data as a prior for the experiment. 
An initial prior distribution is needed and we propose to rely on Jeffreys prior \citep{box_bayesian_1973}. This prior maximizes the distance between the prior and the posterior and is the standard non informative proposal, on the two parameters
$\theta_0$ and $\theta_1$. 
It is defined as:
\begin{equation} \label{prior_jeffreys}
\pi_j^{}\left(\theta_{j}\right)\propto\sqrt{\left\vert I^{F}\left(\theta_{j}\right)\right\vert}=\sqrt{\left\vert -E_{\theta_{j}}\left[\frac{\partial^{2}\log\ell_j\left(\theta_{j}\mid y\right)}{\partial\theta_{j}^{t}\partial\theta_{j}}\right]\right\vert}\quad,\quad j=0,1
\end{equation}
where $I^{F}\left(\theta_{j}\right)$ is Fisher information matrix and is interpreted as the amount of information provided by the observation $y$ on $\theta_{j}$. 
Unfortunately, Jeffreys prior distributions on $\theta_0$ and $\theta_1$ are improper, since the integrals are infinity.
The corresponding posteriors are well-defined but, when an improper prior is used on the parameters, the posterior in the model space is not properly defined \citep{robert_bayesian_2007,gelman_bayesian_2013}. Indeed, the posterior probabilities of both models depend on the unknown normalizing constant.
To compute this posterior, we propose to proceed in three steps using two historical datasets similar to the experiment we want to analyse.
The experiment used in the first step is referred to as the experiment 1 or the first experiment. The experiment used in the second step is referred to as experiment 2 or the second experiment. Finally, the third experiment (also named the experiment 3 or the experiment of interest) used in the last step is the dataset on which the inference is made.
The three steps are:
\begin{enumerate}
\item A first set of historical data, together with the Jeffreys priors on $\theta_0$ and $\theta_1$ are considered to deduce a first proper posterior on $\theta_0$ and $\theta_1$;
\item A second set of historical data is considered with the use the first posterior as prior. We introduce of a Laplace (flat) prior on the two hypotheses, namely $\Pr(\mathcal{M}_0)=\Pr(\mathcal{M}_1)=1/2$, and deduce a new posterior on $\theta_0$ and $\theta_1$ and also on the hypotheses;
\item In the last and third step of our procedure, we deduce the posterior probabilities of the models $\mathcal{M}_0$ and $\mathcal{M}_1$ which are computed using the posterior of the second step as an objective informative prior.
\end{enumerate}
In the next two subsections, the marginal likelihoods for the three step procedure are detailed. The prior and the posterior distributions for each step of our proposal are indexed. 
Let us denote by: 
\begin{itemize}
\item $y_1=(c_1,t_1)$ the first historical set where $c_1$ is an \textit{iid} sample of size $n_{c_1}$ and
$t_1$ an \textit{iid} sample of size $n_{t_1}$;
\item $y_2=(c_2,t_2)$ the second one (of sizes $n_{c_2}$ and $n_{t_2}$);
\item $y_3=(c_3,t_3)$ the data to analyse (of sizes $n_{c_3}$ and $n_{t_3}$);
\item $\overline{c}_{i}$ and $\overline{t}_{i}$ denote the empirical expectations of $c_i$ and $t_i$, for $i=1,2,3$ and $\gamma_{c_{i}}=\sum_{i=1}^{n_{c_{i}}}\left(c_{i}-\overline{c}_{i}\right)^{2}$ and
$\gamma_{t_{i}}=\sum_{i=1}^{n_{t_{i}}}\left(t_{i}-\overline{t}_{i}\right)^{2}$.
\end{itemize}

\subsection{Computation of integrated likelihoods under model $\mathcal{M}_1$}

Let us first consider the case of model $\mathcal{M}_1$ for which all the computations are explicit.
Under model $\mathcal{M}_{1}$, the likelihood of $\theta_1$ for dataset $y_i$ is:
\begin{eqnarray}
\label{vraisemblance_M1}
\ell_1\left(\theta_1|y_i\right) = 
\left(2\pi\sigma_{1,c}^{2}\right)^{\frac{-n_{c_{i}}}{2}} \exp\left\lbrace\frac{-1}{2\sigma_{1,c}^{2}} \left[ n_{c_{i}}\left(\overline{c_{i}}- \mu_{c}\right)^{2} +\gamma_{c_i}\right]\right\rbrace\times \nonumber \\
\left(2\pi\sigma_{1,t}^{2}\right)^{\frac{-n_{t_{i}}}{2}} \exp\left\lbrace\frac{-1}{2\sigma_{1,t}^{2}} \left[ n_{t_{i}}\left(\overline{t_{i}}- \mu_{t}\right)^{2} +\gamma_{t_i}\right]\right\rbrace 
\end{eqnarray}
Regarding the first step of the analysis, the Jeffreys' prior on $\theta_1$ used is:
$$
\pi_1^1\left(\theta_{1}\right)\propto\left(\sigma_{c}^{2}\sigma_{t}^{2}\right)^{-\frac{3}{2}}\,.
$$
The posterior distribution of $\theta_1$ can be then deduced. Indeed,
\begin{eqnarray} \label{eq_posterior_M1_1}
\pi_1^1\left(\theta_1\vert y_1\right) & \propto & \ell_1\left(\theta_1\mid y_1\right)\pi_1^1\left(\theta_1\right) \nonumber \\
                                      & \propto & \pi_1^1\left(\mu_c\vert\sigma_{1,c}^2,y_1\right)\pi_1^1\left(\sigma_{1,c}^2\vert y_1\right) 
                                                  \pi_1^1\left(\mu_t\vert\sigma_{1,t}^2,y_1\right)\pi_1^1\left(\sigma_{1,t}^2\vert y_1\right)
\end{eqnarray}
Equation \eqref{eq_posterior_M1_1} shows posterior independence between parameters of the two groups of the first step.
The conditional posterior distributions of parameters are deduced such as the location parameters $\mu_c$ and $\mu_t$ are Gaussian ($\mathcal{N}$) and the marginal posterior distributions of $\sigma_{1,c}^2$ and $\sigma_{1,t}^2$ belongs to the inverse-gamma family ($\mathcal{I}G$).
Thus, the prior distributions that will be used at the second step of the analysis are:
$$
\mu_c\vert y_1,\sigma_{1,c}^2\sim\mathcal{N}\left(\overline{c_{1}},\frac{\sigma_{1,c}^{2}}{n_{c_1}}\right),\quad
\mu_t\vert y_1,\sigma_{1,t}^2\sim\mathcal{N}\left(\overline{t_{1}},\frac{\sigma_{1,t}^{2}}{n_{t_1}}\right),
$$
$$
\sigma_{1,c}^2\sim\mathcal{I}G\left(\frac{n_{c_1}}{2},\frac{\gamma_{c_1}}{2}\right),\quad
\sigma_{1,t}^2\sim\mathcal{I}G\left(\frac{n_{t_1}}{2},\frac{\gamma_{t_1}}{2}\right)\,.
$$
Using $\pi_1^2(\theta_1)=\pi_1^1(\theta_1|y_1)$, the integrated likelihood of $y_2$ becomes:
\begin{eqnarray}
\Pr(y_2|\mathcal{M}_1) & = & \int \ell_1(\theta_1|y_2)\pi_1^2(\theta_1)d\theta_1 \nonumber\\
                       & = & \int \ell_1(\theta_1|y_2)\pi_1^1(\theta_1|y_1)d\theta_1 \nonumber \\
                       & = & \left(2\pi\right)^{\frac{n_{c_{2}}+n_{t_{2}}}{2}} \left(\frac{n_{c_{1}}n_{t_{1}}}{\left(n_{c_1}+n_{c_2}\right)\left(n_{t_1}+n_{t_2}\right)}\right)^{\frac{1}{2}} 
                       \frac{\left(\frac{\gamma_{t_1}}{2}\right)^{\frac{n_{t_1}}{2}}\left(\frac{\gamma_{c_1}}{2}\right)^{\frac{n_{c_1}}{2}}}
                       {\Gamma\left(\frac{n_{t_1}}{2}\right)\Gamma\left(\frac{n_{c_1}}{2}\right)}\times \nonumber \\
                       &   & \frac{\Gamma\left(\frac{n_{t_1}+n_{t_2}}{2}\right)\Gamma\left(\frac{n_{c_1}+n_{c_2}}{2}\right)}
                             {\left(\frac{\gamma_{t_1}+\gamma_{t_2}}{2}+\frac{n_{t_1}n_{t_2}\left(\overline{t_1}-\overline{t_2}\right)^2}{2\left(n_{t_1}+n_{t_2}\right)}\right)^{\frac{n_{t_1}+n_{t_2}}{2}}
                              \left(\frac{\gamma_{c_1}+\gamma_{c_2}}{2}+\frac{n_{c_1}n_{c_2}\left(\overline{c_1}-\overline{c_2}\right)^2}{2\left(n_{c_1}+n_{c_2}\right)}\right)^{\frac{n_{c_1}+n_{c_2}}{2}}}\,.
\end{eqnarray}
Similarly, using $\pi_1^2(\theta_1)=\pi_1^1(\theta_1|y_1)$, the posterior distribution obtained in the second step is:
\begin{eqnarray} \label{eq_posterior_M1_2}
\pi_1^2\left(\theta_1\vert y_2\right) & \propto & \ell_1\left(\theta_1\mid y_2\right)\pi_1^2\left(\theta_1\right) \nonumber \\
                                      & \propto & \pi_1^2\left(\mu_c\vert\sigma_{1,c}^2,y_2\right)\pi_1^2\left(\sigma_{1,c}^2\vert y_2\right) 
                                                  \pi_1^2\left(\mu_t\vert\sigma_{1,t}^2,y_2\right)\pi_1^2\left(\sigma_{1,t}^2\vert y_2\right).
\end{eqnarray}
Like for the first step, equation \eqref{eq_posterior_M1_2} shows posterior independence between parameters of the two groups for the second step.
Used as prior for the third step, the posterior distribution of each parameter is thus deduced for the second step and corresponds to:
$$
\mu_c\vert y_2,\sigma_{1,c}^2\sim\mathcal{N}\left(\frac{n_{c_1}\overline{c_1}+n_{c_2}\overline{c_2}}{n_{c_1}+n_{c_2}},\frac{\sigma_{1,c}^{2}}{n_{c_{1}}+n_{c_{2}}}\right),\quad
\mu_t\vert y_2,\sigma_{1,t}^2\sim\mathcal{N}\left(\frac{n_{t_1}\overline{t_1}+n_{t_2}\overline{t_2}}{n_{t_1}+n_{t_2}},\frac{\sigma_{1,t}^{2}}{n_{t_{1}}+n_{t_{2}}}\right),
$$
$$
\sigma_{1,c}^2\sim\mathcal{I}G\left(\frac{n_{c_1}+n_{c_2}}{2},\frac{\gamma_{c_1}+\gamma_{c_2}}{2}+
\frac{n_{c_1}n_{c_2}\left(\overline{c_1}-\overline{c_2}\right)^2}{2\left(n_{c_1}+n_{c_2}\right)}\right),
$$
$$
\sigma_{1,t}^2\sim\mathcal{I}G\left(\frac{n_{t_1}+n_{t_2}}{2},\frac{\gamma_{t_1}+\gamma_{t_2}}{2}+
\frac{n_{t_1}n_{t_2}\left(\overline{t_1}-\overline{t_2}\right)^2}{2\left(n_{t_1}+n_{t_2}\right)}\right).
$$
Thus, using $\pi_1^3(\theta_1)=\pi_1^2(\theta_1|y_2)$, the integrated likelihood of $y_3$ becomes:
\begin{eqnarray}
\Pr(y_3|\mathcal{M}_1) & = & \int \ell_1(\theta_1|y_3)\pi_1^3(\theta_1)d\theta_1 \nonumber\\
                       & = & \int \ell_1(\theta_1|y_3)\pi_1^2(\theta_1|y_2)d\theta_1 \nonumber \\
                       & = & \left(2\pi\right)^{\frac{n_{c_3}+n_{t_3}}{2}}\left(\frac{\left(n_{c_1}+n_{c_2}\right)\left(n_{t_1}+n_{t_2}\right)}
                             {\left(n_{c_1}+n_{c_2}+n_{c_3}\right)\left(n_{t_1}+n_{t_2}+n_{t_3}\right)}\right)^{\frac{1}{2}} 
                             \frac{\left(\frac{\gamma_{t_1}+\gamma_{t_2}}{2}+\frac{n_{t_1}n_{t_2}\left(\overline{t_1}-\overline{t_2}\right)^2}{2\left(n_{t_1}+n_{t_2}\right)}\right)^{\frac{n_{t_1}+n_{t_2}}{2}}}
                             {\Gamma\left(\frac{n_{t_1}+n_{t_2}}{2}\right)} \nonumber \\
                       &   & \times\frac{\left(\frac{\gamma_{c_1}+\gamma_{c_2}}{2}+\frac{n_{c_1}n_{c_2}\left(\overline{c_1}-\overline{c_2}\right)^2}{2\left(n_{c_1}+n_{c_2}\right)}\right)^{\frac{n_{c_1}+n_{c_2}}{2}}}
                             {\Gamma\left(\frac{n_{c_1}+n_{c_2}}{2}\right)} \frac{\Gamma\left(\frac{n_{t_1}+n_{t_2}+n_{t_3}}{2}\right)\Gamma\left(\frac{n_{c_1}+n_{c_2}+n_{c_3}}{2}\right)}
                             {\beta_t^{\frac{n_{t_1}+n_{t_2}+n_{t_3}}{2}}\beta_c^{\frac{n_{c_1}+n_{c_2}+n_{c_3}}{2}}}    \,  ,           
\end{eqnarray}
where
$$
\beta_t=\frac{\gamma_{t_1}+\gamma_{t_2}+\gamma_{t_3}}{2}+
\frac{n_{t_1}n_{t_2}\left(\overline{t_1}-\overline{t_2}\right)^2}{2\left(n_{t_1}+n_{t_2}\right)}
+\frac{\left(n_{t_1}+n_{t_2}\right)n_{t_3}}{2\left(n_{t_1}+n_{t_2}+n_{t_3}\right)}
\left(\frac{n_{t_1}\overline{t_1}+n_{t_2}\overline{t_2}}{n_{t_1}+n_{t_2}}-\overline{t_3}\right)^{2},
$$
$$
\beta_c=\frac{\gamma_{c_1}+\gamma_{c_2}+\gamma_{c_3}}{2}+
\frac{n_{c_1}n_{c_2}\left(\overline{c_1}-\overline{c_2}\right)^2}{2\left(n_{c_1}+n_{c_2}\right)}
+\frac{\left(n_{c_1}+n_{c_2}\right)n_{c_3}}{2\left(n_{c_1}+n_{c_2}+n_{c_3}\right)}
\left(\frac{n_{c_1}\overline{c_1}+n_{c_2}\overline{c_2}}{n_{c_1}+n_{c_2}}-\overline{c_3}\right)^{2}\,.
$$

\subsection{Calculation of integrated likelihoods under model $\mathcal{M}_0$}

Under model $\mathcal{M}_{0}$, the likelihood of $\theta_0$ for dataset $y_i$ is:
\begin{eqnarray}
\ell_0\left(\theta_0| y_i\right) & = & \left( 2\pi\sigma_{0,c}^{2}\right)^{\frac{-n_{c_i}}{2}} \left( 2\pi\sigma_{0,t}^{2}\right)^{\frac{-n_{t_i}}{2}}\nonumber \\
& & \times\exp\left\lbrace-\frac{\left( n_{c_i}\sigma_{0,t}^{2} +n_{t_i}\sigma_{0,c}^{2}\right)}{2\sigma_{0,c}^{2}\sigma_{0,t}^{2}}  \left(\mu -\frac{n_{c_i}\sigma_{0,t}^{2}\overline{c_i} +n_{t_i}\sigma_{0,c}^{2}\overline{t_i}}{n_{c_i}\sigma_{0,t}^{2} +n_{t_i}\sigma_{0,c}^{2}}\right)^{2}\right\rbrace\nonumber\\
& & \times\exp\left\lbrace-\frac{\gamma_{t_i}}{2\sigma_{0,t}^{2}}-\frac{\gamma_{c_i}}{2\sigma_{0,c}^{2}} -\frac{n_{c_i} n_{t_i}\left(\overline{c_i} -\overline{t_i}\right)^{2}}{2\left(n_{c_i}\sigma_{0,t}^{2} +n_{t_i}\sigma_{0,c}^{2}\right)} \right\rbrace . \nonumber
\end{eqnarray}
The Jeffreys prior on $\theta_0$ used in the first step of the analysis is:
$$
\pi_0^1\left(\theta_{0}\right)\propto\left(\frac{n_{c}\sigma_{0,t}^{2}+n_{t}\sigma_{0,c}^{2}}{\sigma_{0,c}^{6}\sigma_{0,t}^{6}}\right)^{\frac{1}{2}}\,.
$$
For each step in the analysis, 
the posterior distribution of $\mu$ is a Gaussian distribution with expectations and variances described in Table \ref{tab_posteriorM0_mu}.
However, due to the dependence between $\sigma_{0,c}^2$ and $\sigma_{0,t}^2$, 
the posterior distribution of the couple $\left(\sigma_{0,c}^2,\sigma_{0,t}^2\right)$ is not explicit.
For instance, for the first step of the analysis, the posterior density of $(\sigma_{0,t}^2,\sigma_{0,t}^2)$ is:
\begin{equation}
\pi_0^1\left(\sigma_{0,c}^2,\sigma_{0,t}^2|y_1\right)\propto
\left(\sigma_{0,c}^{2}\right)^{\frac{-n_{c_{1}}}{2}-1} \left(\sigma_{0,t}^{2}\right)^{\frac{-n_{t_{1}}}{2}-1}
\exp\left\lbrace -\frac{\gamma_{c_1}}{2\sigma_{0,c}^{2}}-\frac{\gamma_{t_1}}{2\sigma_{0,t}^{2}}-
\frac{n_{c_1}n_{t_1}\left(\overline{c_1} -\overline{t_1}\right)^{2}}{2\left(n_{c_1}\sigma_{0,t}^{2} +n_{t_1}\sigma_{0,c}^{2}\right)} \right\rbrace\,.
\label{post_var}
\end{equation}

\begin{table}[!th]
\caption{\label{tab_posteriorM0_mu}
Under model $\mathcal{M}_0$, posterior expectation and variance of $\mu$ conditional on $\sigma_{0,c}^2$ and $\sigma_{0,t}^2$ for the three steps of the analysis.}
\begin{center}
\scalebox{0.7}{
\begin{tabular}{|c||c|c|}
\hline
\multicolumn{1}{|c||}{$\mu|\sigma_{0,c}^2,\sigma_{0,t}^2,y_i$} & \multicolumn{1}{c|}{$\E\left[\mu|\sigma_{0,c}^2,\sigma_{0,t}^2,y_i\right]$} & 
\multicolumn{1}{c|}{$\V\left[\mu|\sigma_{0,c}^2,\sigma_{0,t}^2,y_i\right]$}\\
\hline
\hline
$\underset{\text{(conditioned on $y_1$, first historical data)}}{\text{Step 1}} $ & $\ds \frac{n_{c_{1}}\sigma_{0,t}^{2}\overline{c_{1}}+n_{t_{1}}\sigma_{0,c}^{2}\overline{t_{1}}}{n_{c_{1}}\sigma_{0,t}^{2}+n_{t_{1}}\sigma_{0,c}^{2}}$ & $\ds \frac{\sigma_{0,c}^{2}\sigma_{0,t}^{2}}{n_{c_{1}}\sigma_{0,t}^{2}+n_{t_{1}}\sigma_{0,c}^{2}}$ \\
\hline
$\underset{\text{(conditioned on $y_2$, second historical data)}}{\text{Step 2}} $  & $\ds \frac{\sigma_{0,t}^{2}\left(n_{c_{1}}\overline{c_{1}}+n_{c_{2}}\overline{c}_{2}\right) +\sigma_{0,c}^{2}\left(n_{t_{1}}\overline{t}_{1}+n_{t_{2}}\overline{t}_2 \right)} {\sigma_{0,t}^{2}\left(n_{c_{1}}+n_{c_{2}}\right)+\sigma_{0,c}^{2}\left(n_{t_{1}}+n_{t_{2}}\right)}$ & $\ds \frac{\sigma_{0,c}^{2}\sigma_{0,t}^{2}}{\sigma_{0,t}^{2}\left(n_{c_{1}}+n_{c_{2}}\right)+\sigma_{0,c}^{2}\left(n_{t_{1}}+n_{t_{2}}\right)}$ \\
\hline
$\underset{\text{(conditioned on $y_3$, interested experiment)}}{\text{Step 3}} $  & $\ds \frac{\sigma_{0,t}^{2}\left(n_{c_{1}}\overline{c_{1}}+n_{c_{2}}\overline{c}_{2}+n_{c_{3}}\overline{c}_{3}\right) +\sigma_{0,c}^{2}\left(n_{t_{1}}\overline{t}_{1}+n_{t_{2}}\overline{t}_2 +n_{t_{3}}\overline{t}_3\right)} {\sigma_{0,t}^{2}\left(n_{c_{1}}+n_{c_{2}}+n_{c_{3}}\right)+\sigma_{0,c}^{2}\left(n_{t_{1}}+n_{t_{2}}+n_{c_{2}} \right)} $ & $\ds \frac{\sigma_{0,c}^{2}\sigma_{0,t}^{2}}{\sigma_{0,t}^{2}\left(n_{c_{1}}+n_{c_{2}}+n_{c_{3}}\right)+\sigma_{0,c}^{2}\left(n_{t_{1}}+n_{t_{2}}+n_{c_{2}} \right)} $ \\
\hline
\end{tabular}
}
\end{center}
\end{table}

A Markov Chain Monte Carlo (MCMC) algorithm through the software WinBUGS \citep{cowles_review_2004,kery_introduction_2010} is used to obtain
a posterior sample from (\ref{post_var}). 
Then, as prior for the second and third steps of the analysis, a product of Inverse-Gamma distributions with parameters calibrated using the posterior sample given by WinBUGS is considered. According to the same process for the steps of the analysis $i=1,2$, this corresponds to:
\begin{eqnarray*}
\underbrace{\pi_0^i\left(\sigma_{0,c}^{2},\sigma_{0,t}^{2}\vert y_i\right)}_{unknown} \xrightarrow[\text{Independence hypothesis}]{\text{MCMC methods}} \begin{array}{c} 
\pi_0^i\left(\sigma_{0,c}^{2}\vert y_i\right)=\pi_0^i\left(\sigma_{0,c}^2\vert\widehat{\alpha}_{c_i},\widehat{\beta}_{c_i}\right)\sim\mathcal{IG}\left(\widehat{\alpha}_{c_i},\widehat{\beta}_{c_i}\right)\\
\pi_0^i\left(\sigma_{0,t}^{2}\vert y_i\right)=\pi_0^2\left(\sigma_{0,t}^2\vert\widehat{\alpha}_{t_i},\widehat{\beta}_{t_i}\right)\sim\mathcal{IG}\left(\widehat{\alpha}_{t_i},\widehat{\beta}_{t_i}\right)
\end{array},
\end{eqnarray*}
where the parameters $\widehat{\alpha}$ and $\widehat{\beta}$ are estimated by the maximum likelihood method
on the MCMC paths. Note that, for the second step of theanalysis, we have the posterior density of $(\sigma_{0,c}^2,\sigma_{0,t}^2)$, target of the second MCMC scheme, such as:
{\footnotesize \begin{eqnarray*}
\pi_0^2\left(\sigma_{0,c}^2,\sigma_{0,t}^2|y_2\right) & \propto & \left(\sigma_{0,c}^{2}\right)^{\frac{-n_{c_{2}}}{2}-\widehat{\alpha}_{c_1}-1} \left(\sigma_{0,t}^{2}\right)^{\frac{-n_{t_{1}}}{2}-\widehat{\alpha}_{t_1}-1} 
\exp\left\lbrace -\frac{\widehat{\beta}_{c_1}}{\sigma_{0,c}^{2}} \right\rbrace \exp\left\lbrace -\frac{\widehat{\beta}_{t_1}}{\sigma_{0,t}^{2}} \right\rbrace   \left(\frac{\sigma_{0,t}^{2}n_{c_{1}}+\sigma_{0,c}^{2}n_{t_{1}}}{\sigma_{0,t}^{2}\left(n_{c_{1}}+n_{c_{2}}\right)+\sigma_{0,c}^{2}\left(n_{t_{1}}+n_{t_{2}}\right)}\right)^{\frac{1}{2}} \\
& & \times \exp\left\lbrace -\frac{\left(n_{c_{1}}\sigma_{0,t}^{2}+n_{t_{1}}\sigma_{0,c}^{2}\right) \left(n_{c_{2}}\sigma_{0,t}^{2}+n_{t_{2}}\sigma_{0,c}^{2}\right)}{2\sigma_{0,c}^{2}\sigma_{0,t}^{2} \left[\sigma_{0,t}^{2}\left(n_{c_{1}}+n_{c_{2}}\right)+\sigma_{0,c}^{2}\left(n_{t_{1}}+n_{t_{2}}\right)\right]} \left( \frac{n_{c_{2}}\sigma_{0,t}^{2}\overline{c_{2}}+n_{t_{2}}\sigma_{0,c}^{2}\overline{t_{2}}}{n_{c_{2}}\sigma_{0,t}^{2}+n_{t_{2}}\sigma_{0,c}^{2}}-\E\left[\mu|\sigma_{0,c}^2,\sigma_{0,t}^2,y_1\right] \right)^{2} \right\rbrace \\
& & \times \exp\left\lbrace -\frac{\gamma_{c_2}}{2\sigma_{0,c}^{2}} -\frac{\gamma_{t_2}}{2\sigma_{0,t}^{2}} -\frac{n_{c_{2}}n_{t_{2}}\left(\overline{c_{2}} -\overline{t_{2}}\right)^{2}}{2\left(n_{c_{2}}\sigma_{0,t}^{2} +n_{t_{2}}\sigma_{0,c}^{2}\right)} \right\rbrace \,.
\end{eqnarray*}}
For steps 2 and 3, the prior distributions are:
\begin{equation*}
\pi_0^2\left(\theta_{0}\right)=\pi_0^1\left(\mu\vert\sigma_c^2,\sigma_t^2,y_1\right)
\pi_0^1\left(\sigma_c^2\vert\widehat{\alpha}_{c_1},\widehat{\beta}_{c_1}\right)
\pi_0^1\left(\sigma_t^2\vert\widehat{\alpha}_{t_1},\widehat{\beta}_{t_1}\right),
\end{equation*}
\begin{equation}
\pi_0^3\left(\theta_{0}\right)=\pi_0^2\left(\mu\vert\sigma_c^2,\sigma_t^2,y_2\right)
\pi_0^2\left(\sigma_c^2\vert\widehat{\alpha}_{c_2},\widehat{\beta}_{c_2}\right)
\pi_0^2\left(\sigma_t^2\vert\widehat{\alpha}_{t_2},\widehat{\beta}_{t_2}\right). \label{eq_prior3_M0}
\end{equation}
Thus, the integrated likelihoods are:
\begin{eqnarray*}
\Pr\left(y_{2}\vert\mathcal{M}_{0}\right) & = & \int\ell_0\left(\theta_0\mid  y_2\right)\pi_0^2\left(\theta_0\right) d\theta_{0} \\
										  & = & \int\ell_0\left(\theta_0\mid  y_2\right)\pi_0^1\left(\theta_0|y_1\right) d\theta_{0} \\
                                          & = & \int\ell_0\left(\theta_0\mid  y_2\right)\pi_0^1\left(\mu\vert\sigma_{0,c}^2,\sigma_{0,t}^2,y_1\right)
                                                \pi_0^1\left(\sigma_{0,c}^2\vert\widehat{\alpha}_{c_1},\widehat{\beta}_{c_1}\right)
                                                \pi_0^1\left(\sigma_{0,t}^2\vert\widehat{\alpha}_{t_1},\widehat{\beta}_{t_1}\right)d\mu d\sigma_{0,c}^{2} d\sigma_{0,t}^{2} 
\end{eqnarray*}
\begin{eqnarray*}
\Pr\left(y_{3}\vert\mathcal{M}_{0}\right) & = & \int\ell_0\left(\theta_0\mid  y_3\right)\pi_0^3\left(\theta_0\right) d\theta_{0} \\
										  & = & \int\ell_0\left(\theta_0\mid  y_3\right)\pi_0^2\left(\theta_0|y_2\right) d\theta_{0} \\
                                          & = & \int\ell_0\left(\theta_0\mid  y_3\right)\pi_0^2\left(\mu\vert\sigma_{0,c}^2,\sigma_{0,t}^2,y_2\right)
                                                \pi_0^2\left(\sigma_{0,c}^2\vert\widehat{\alpha}_{c_1},\widehat{\beta}_{c_1}\right)
                                                \pi_0^2\left(\sigma_{0,t}^2\vert\widehat{\alpha}_{t_1},\widehat{\beta}_{t_1}\right)d\mu d\sigma_{0,c}^{2} d\sigma_{0,t}^{2} 
\end{eqnarray*}
The integrations over $\mu$ are explicit, in contrast to the variance parameters. Numerical integration is thus used to solve this difficulty through the R package \textit{cubature} and the \textit{adaptInegrate} function which implements an adaptive quadrature method \citep{cubature}.

\section{Simulation study}

In this section, the relative efficiency of the proposed method is studied and compared with classical approaches concerning the Behrens-Fisher problem using simulation.
First, we assess the frequentist properties of our Bayesian approach as recommended by the FDA \citep{food_and_drug_administration_guidance_2010}. 
Our methodology is compared with two currently used approaches: Student test with Satterthwaite correction applied only on the experiment of interest and Student test with Satterthwaite correction applied on the three experiments pooled.
Using the frequestist approaches, the null hypothesis is rejected if the \textit{p-value} is less than $5\%$.
Then, our method's behavior is also assessed in the similar cases when a weight is affected to historical data.

The output of our Bayesian method is the posterior probability $Pr\left(\mathcal{M}_{1}\vert y\right)$ which is directly interpretable and more intuitive than a \textit{p-value}.
Given this probability, the experimenter has to discriminate between the two hypotheses. Thus, a decision rule has to be defined according to a chosen threshold.
In our case, we accept that the expectations are different if the posterior probability is greater than a certain threshold $p$ such as $Pr\left(\mathcal{M}_{1}\vert y\right)>p$. 
The first intuitive threshold is $p=0.5$. 
Indeed, in the case of two models, if the posterior probability of one model is greater than $0.5$, this model is more probable than the other one, given the data.
We also propose another arbitrary threshold which is more conservative for $\mathcal{M}_0$: $p=0.8$. 
If the posterior probability of the model is greater than $0.8$, we assume the probability is enough to decide with a minimal risk of concluding wrongly to a difference.
Finally, a threshold (tailor-made) was also calibrated to control the type I error rate given the studied situations.
From the simulations under the model $\mathcal{M}_0$, the latter threshold corresponded to the posterior probability $Pr\left(\mathcal{M}_{1}\vert y\right)$ associated with the $5\%$ rank.

We simulated three datasets (experiments) using Gaussian distributions in the context of the Complete Freund's Adjuvant (CFA) protocol, detailed in the real data analysis section.
The values of the parameters in the simulations are chosen according to prior knowledge gained from experimental results in the CFA protocol. 
The small sample sizes $\left(n_{c_i}\right)_{i=1,2,3}=10$ and $\left(n_{t_i}\right)_{i=1,2,3}=10$ are considered.
From all the simulated data, we defined $\mu_c=2.94$ and $\mu_t=\mu_c+\mu_c\delta$ where $\delta$ is the percentage of effect.
Concerning the standard deviation of the control group, $\left(\sigma_{c_i}\right)_{i=1,2,3}=0.6$ whatever the situation. Regarding the treated group standard deviations $\left(\sigma_{t_i}\right)_{i=1,2,3}$, we considered the different values $0.6$, $1.5$ and $3$, respectively associated with the minimum, median and maximum values deduced from historical data.
For all simulation studies, each experiment is simulated a thousand times ($N=1000$).

\subsection{Comparison between the Bayesian and frequentist approaches}

The first situation presented in Table \ref{tab_simu1} is the ideal case where the three experiments are similar. 
The expectations of the treated group were deduced from several effects relative the control group: $0\%$ for the verification of the type I error, then $30\%$, $40\%$ and $50\%$. 
We began to study a difference of $30\%$ because this was considered to be the minimum meaningful effect in the biological experts' opinion.
Both historical datasets were similar with the same standard deviation for control groups $\left(\sigma_{c_i}=0.6\right)_{i=1,2}$ and treated groups $\left(\sigma_{t_i}=1.5\right)_{i=1,2}$.
Concerning the experiment of interest (experiment 3), the different standard deviations were studied for the treated group.
When there was no difference between the expectations of treated and control groups, both frequentist approaches controlled the type I error rate: the difference between the two expectations was detected in maximum $5\%$ of cases. 
Regarding the Bayesian approach, the choice of the threshold to ensure a sufficient posterior probability of the model $\mathcal{M}_1$ was important. Indeed, a threshold of $0.5$ did not seem strict enough because the type I error rate was not controlled: we concluded to a difference between the two expectations in $10\%$ to $20\%$ of cases while it was not existant.
However, a threshold of $0.8$ allowed a better control of the type I error rate than for the frequentist approaches, and even if it is conservative for the $\mathcal{M}_0$.
Whatever the considered standard deviation for the treated group in the third experiment, the calibrated threshold verifying this error to $5\%$ did not exceed $0.7$.

When the standard deviation of experiment 3 treated group varied from those of experiments 1 and 2 (Table \ref{tab_simu1}), the power of the three approaches was affected. Indeed, it decreased when treated group variability in step 3 increased. 
In all cases, the Bayesian approach was better than the Student test performed only on the experiment 3 for each variation. 
The frequentist approach where data from the three experiments were pooled was the most powerful.
Indeed, this approach was very robust in this case because to pool data from Gaussian distribution with the same expectation was like considering one sample with an expectation and a variance from the three experiment's variances. 
Thus, the Student test was the best approach because it was in the ideal case: two Gaussian samples with the size of $30$ unit.
However, the experiment of interest is the last and the two first are taken into account from an informative point of view.  
We want to conclude only on the last experiment and not the all of them. 

\begin{table}[h!]
\caption{Efficiency (in percentage) of methods to detect a difference between the two sample expectations when $\sigma_{t_3}$ changes given a difference of $\delta$ percent. 
The thresholds $p$ to respect the type I error rate to $5\%$ for the values of $\sigma_{t_{3}}$ min, median and max are respectively 0.672, 0.698 and 0.660.}
\centering
\scalebox{0.95}{
\begin{tabular}{|l||c|c|c||c|c|c||c|c|c||c|c|c|} 
\hline
\multicolumn{1}{|c||}{$\delta$} & \multicolumn{3}{c||}{$0\%$} & \multicolumn{3}{c||}{$30\%$} & \multicolumn{3}{c||}{$40\%$} & \multicolumn{3}{c|}{$50\%$}\\ 
\hline
\multicolumn{1}{|c||}{$\sigma_{t_{3}}$}  & \multicolumn{1}{c|}{min}  & \multicolumn{1}{c|}{med} & \multicolumn{1}{c||}{max} & \multicolumn{1}{c|}{min} & \multicolumn{1}{c|}{med} & \multicolumn{1}{c||}{max}  & \multicolumn{1}{c|}{min}  & \multicolumn{1}{c|}{med} & \multicolumn{1}{c||}{max} & \multicolumn{1}{c|}{min}  & \multicolumn{1}{c|}{med} & \multicolumn{1}{c|}{max}  \\ 
\hline
T.test & 4.3 & 5.0 & 5.0 & 87.1 & 34.5 & 14.4 & 98.5 & 59.6 & 18.3 & 100 & 76.7 & 29.8 \\
T.test (pool)  & 5.4 & 4.3 & 5.3 & 91.4 & 81.2 & 58.6 & 99.3 & 97.2 & 79.7 & 99.9 & 99.8 & 94.5 \\
$\Pr\left(\mathcal{M}_{1}\vert y\right)>0.5$ & 16.4 & 17.1 & 12.7 & 94.1 & 82.8 & 53.0 & 98.8 & 94.4 & 63.8 & 100 & 98.3 & 74.3\\
$\Pr\left(\mathcal{M}_{1}\vert y\right)>p$ & 4.9  & 5  & 5  & 88.3 & 68.8 & 40.3 & 97.3 & 89.2 & 53.5 & 100 & 96.1 & 65.2\\
$\Pr\left(\mathcal{M}_{1}\vert y\right)>0.8$ & 1.3 & 1.7 & 2.3 & 76.7 & 59.1 & 29.1 & 95.5 & 84.8 & 42.0 & 100 & 94.2 & 56.4\\
\hline
\end{tabular}
}
\label{tab_simu1}
\end{table}

To pool datasets can influence wrongly the statistical test. 
The Bayesian approach uses information from the two previous experiments to calculate the posterior probability of each model about the third experiment.
Table \ref{tab_simu_loc1} presents the results from different situations which illustrate the historical data influence when the location parameters $\left(\mu_{t_i}\right)_{i=1,2,3}$ are different between experiments $i$. 
For each situation, the difference between the control group expectation and the treated group expectation were distinct by experiments. 
The Student test performed only on the experiment of interest was not influenced by the historical datasets because it did not take them into account.
However, the Bayesian approach and the application of the Student test on the pooled data were affected by historical data information, and concluded sometimes wrongly on the difference concerning the experiment of interest.
In situations 6 to 9, there was no distinction between treated and control group concerning the third experiment, in contradiction with historical experiments. 
In these cases, the approach with pooled data gave the same weight between the three experiments concluding wrongly with a high percentage. 
Thereby bringing no difference between the three experiments which was not our aim in this work.
The consequence of an insignificant effect for the first two experiments (first two columns) was the decrease of the method's power, compared to the previous results where the three experiments were similar. 
The Bayesian approach made a difference between the data of each experiment, giving the conclusion on the last experiment. 
We noticed a natural weighting for our method where the conclusion on the last experiment was less influenced by the first historical experiment than the second one.
Indeed, when the difference between the two groups changed only on the experiment considered in the second step (situation 8 versus 9), the percentage to conclude wrongly on the difference varied more strongly than when this change was performed only on experiment of the first step (situation 6 versus 8).
As mentioned before, considering a threshold of $0.5$ did not allow to respect the good frequentist properties. Indeed, the type I error rate was not controlled to $5.4\%$. 
The threshold $p=0.8$ seemed more relevant to ensure the decision. 
The threshold $p=0.698$ was also relevant to ensure the decision as the type I error rate was relatively controlled while the power was higher than when using $p=0.8$. 
The historical data has to be chosen appropriately so that the information objectively helps the decision.
Thus, the Bayesian approach is to be the golden mean of the two frequentist approaches using the Student test.
Indeed, the information from historical experiments impacted on the posterior probability with a balanced influence where the experiment used in the second step was more informative than the first one associated with the first step.
The Bayesian method takes into account the information provided by the historical data. It is important to pay attention to the experiments set \textit{a priori} as they affect the results.
The Student test on pooled data is not suitable for our aim, it is not taken into account in the next studies.

\begin{table}[!h]
\caption{Efficiency (in percentage) of methods to detect a difference between the two sample expectations when the location parameter $\left(\mu_{t_i}\right)_{i=1,2,3}$ changes from $\left(\mu_{c_i}\right)_{i=1,2,3}=2.94$ between experiments following the percentage $\delta$. 
The variability of three experiments being equivalent, the threshold associated with the posterior probability is $p=0.698$.}
\centering
\begin{tabular}{|l|c|c|c|c|c|c|c|c|c|c|} 
\hline
\multicolumn{1}{|c|}{\textbf{Situations}} & \textbf{1} & \textbf{2} & \textbf{3} & \textbf{4} & \textbf{5} & \textbf{6} & \textbf{7} & \textbf{8} & \textbf{9} \\
\hline
\multicolumn{1}{|c|}{$\delta$ of Experiment 1} & $10\%$ & $10\%$ & $30\%$ & $50\%$ & $50\%$ & $50\%$ &$50\%$ & $10\%$ & $10\%$ \\
\multicolumn{1}{|c|}{$\delta$ of Experiment 2} & $10\%$ & $10\%$ & $30\%$ & $30\%$ & $10\%$ & $10\%$ &$30\%$ & $10\%$ & $20\%$ \\
\multicolumn{1}{|c|}{$\delta$ of Experiment 3} & $40\%$ & $50\%$ & $10\%$ & $10\%$ & $10\%$ & $0\%$ & $0\%$ & $0\%$ & $0\%$ \\
\hline
\hline
 T.test & 53.7 &  74.9 & 8.3 & 8.2 & 8.9 & 5.4 & 5.5 & 4.3 & 6.1\\
 T.test (pool)  &  44.1 & 58.1 & 61 & 80.7 & 56.6 & 44.0 & 66.2 &  9.6 & 16.2 \\
 $\Pr\left(\mathcal{M}_{1}\vert y\right)>0.5$ & 67.4 & 78.5 & 57.3 & 45.9 & 21.8 & 13.7  & 30.8 & 20.7 & 31.9 \\
$\Pr\left(\mathcal{M}_{1}\vert y\right)>p$ & 44.7 & 59.4 & 38.7 & 31.9 & 13.1 & 6.1  & 19.9 & 6.8 & 13.1\\
$\Pr\left(\mathcal{M}_{1}\vert y\right)>0.8$ & 32.4 & 45.6 & 29.0 & 24.1 & 8.9 & 3.4 & 15.2 & 3.5 & 7.9\\
\hline
\end{tabular}
\label{tab_simu_loc1}
\end{table} 

Table \ref{tab_simu_var12} shows the influence of the variability of two prior experiments on the Bayesian method.
Using a median variability for the prior experiments, the results are similar to those in Table \ref{tab_simu1} where $\sigma_{t_3}=1.5$ (median value).
Also, the variability modification of these experiments did not change Student test power because it was applied only on the third experiment.
The power of the Bayesian method seemed to increase when the variability increases for the experiment 1. In this case, the prior of the first step was more vague and allowed to calibrate an objective prior more efficiently in the second step.
However, the Bayesian method power decreased when the variability increased for the second experiment. Indeed, the weight of the second experiment was more important than the first, and a prior less diffused improved the method power.
Thereby, the variability of the second experiment impacted the power more importantly than the first experiment variability.

\begin{table}[!h]
\caption{Efficiency (in percentage) of methods to detect a difference between the two sample expectations given a difference of $\delta$ percent, when $(\sigma_{t_i})_{i=1,2}$ changes, 
and $\sigma_{t_3}=\left(\sigma_{t_k}\right)_{k\neq i}=1.5$. 
}
\centering
\scalebox{1}{
\begin{tabular}{|c|l||c|c||c|c||c|c||c|c|} 
\hline
\multicolumn{2}{|c||}{$\delta$} & \multicolumn{2}{c||}{$0\%$} & \multicolumn{2}{c||}{$30\%$} & \multicolumn{2}{c||}{$40\%$} & \multicolumn{2}{c|}{$50\%$} \\ 
\hline
\multicolumn{2}{|c||}{$\sigma_{t_{i}}$}  & \multicolumn{1}{c|}{min}  & \multicolumn{1}{c||}{max} & \multicolumn{1}{c|}{min}  & \multicolumn{1}{c||}{max} & \multicolumn{1}{c|}{min}  & \multicolumn{1}{c||}{max}  & \multicolumn{1}{c|}{min}  & \multicolumn{1}{c|}{max} \\ 
\hline
\hline
$i=1$& $\Pr\left(\mathcal{M}_1\vert y\right)>0.5$ & 8.4 & 13.4 & 57.2 & 79.0 & 75.0 & 93.5 &  87.3 & 97.9 \\
& $\Pr\left(\mathcal{M}_1\vert y\right)>p$ & 5 ($p$=0.59) & 5 ($p$=0.618) & 53.4 & 66.6 & 71.9 & 87.9 & 84.9 & 95.2 \\
 &  $\Pr\left(\mathcal{M}_1\vert y\right)>0.8$ & 1.8 & 1.0 & 40.1 & 43.5 & 61.7 & 71.1 & 79.0 & 88.8 \\
\hline
$i=2$ &  $\Pr\left(\mathcal{M}_1\vert y\right)>0.5$ & 13.6 & 15.2 & 93.2 & 53.5 & 98.7 & 64.2 & 99.9 & 74.5  \\
&  $\Pr\left(\mathcal{M}_1\vert y\right)>p$ & 5 ($p$=0.67) & 5 ($p$=0.703) & 86.8 & 36.8 & 97.6 & 47.7 & 99.5 & 63.5 \\
&  $\Pr\left(\mathcal{M}_1\vert y\right)>0.8$ & 1.2 & 2.0 & 77.4 & 28.1 & 94.7 & 41.3 & 99.1 & 56.2  \\
\hline
\end{tabular}
}
\label{tab_simu_var12}
\end{table}

Our method was sensitive to prior experiments. It is important to choose similar experiments where the variability and treatment effect are close. This proximity between experiments is all the more important between the second experiment and the experiment of interest.

\subsection{Influence of weighting for historical experiments}

Although, the three step procedure includes a natural weighting through the updating step of the posterior probability calculation in the model space, an additional weighting could be applied to decrease the prior information influence. 
The aim is to downweight the historical data to some degree \citep{ibrahim_power_2000},
by for example reducing the sample size using a multiplication factor between 0 and 1 given the same empirical expectations and variances.
We characterised this weighting (in percent) with the couple $\left(w_1,w_2\right)$ being the two weights respectively associated with experiment 1 and 2. 
They allowed to decrease the sample size of these prior experiments in the computation and caused more vague prior distributions on the parameters. 
In prior distributions for the third step ($\pi^3_0\left(\theta_0\right)$ for $\mathcal{M}_1$ with the equation \eqref{eq_posterior_M1_2}, and $\pi^3_0\left(\theta_0\right)$ for $\mathcal{M}_0$ with the equation \eqref{eq_prior3_M0}), considering the quantities $\left(w_1n_{c_1},w_1n_{t_1},w_2n_{c_2},w_2n_{t_2}\right)$ instead of $\left(n_{c_1},n_{t_1},n_{c_2},n_{t_2}\right)$ impacted the deviation parameters of the prior distributions (Gaussian and inverse gamma distributions) and caused diffused priors.  
The previous simulation results showed a natural weighting between the two historical experiments: the experiment used in the second step of the Bayesian method influenced more than the one used in the first step. 
Thus, this natural weighting is pertinent and suitable. In this subsection, we considered an additional weighting with the same weights for both historical experiments ($i=1,2$) such as $w=w_1=w_2$.
The previous simulation studies were performed once more with the two following weights: $w=1/3$ and $w=1/2$.

Table \ref{tab_simu1_w} is equivalent to Table \ref{tab_simu1} with the weighting. 
The thresholds calibrated to control the type I error rate to $5\%$ increased with the incertitude. Indeed, when the three experiments were similar, these thresholds were 0.8 and 0.9 respectively for the weights  $w=1/2$ and $w=1/3$.
Since the percentage to conclude wrongly on the different expectations increased with the incertitude, the lower the weight was, the higher the threshold to be considered.
Therefore, regarding $\delta>0$, the results were close with and without weights $w$. With weight, the percentage to detect a difference between expectation seemed to increase for a higher value of $\sigma_{t_3}$ and seemed to decrease for a lower value of $\sigma_{t_3}$ in comparison with results of Table \ref{tab_simu1}.

\begin{table}[!h]
\caption{Efficiency (in percentage) of methods to detect a difference between the two sample expectations when $\sigma_{t_3}$ changes, given a difference of $\delta$ percent and the weighting $w$. The thresholds $p$ to respect the type error to $5\%$ for the values of $\sigma_{t_{3}}$ min, median and max are respectively 0.82, 0.89.6 and 0.73 for $w=\frac{1}{3}$ and respectively 0.753, 0.798, 0.669 for $w=\frac{1}{2}$.}
\centering
\scalebox{0.85}{
\begin{tabular}{|l||c|c|c||c|c|c||c|c|c||c|c|c|} 
\hline
\multicolumn{1}{|c||}{$\delta$} & \multicolumn{3}{c||}{$0\%$} & \multicolumn{3}{c||}{$30\%$} & \multicolumn{3}{c||}{$40\%$} & \multicolumn{3}{c|}{$50\%$}\\ 
\hline
\multicolumn{1}{|c||}{$\sigma_{t_{3}}$}  & \multicolumn{1}{c|}{min}  & \multicolumn{1}{c|}{median} & \multicolumn{1}{c||}{max} & \multicolumn{1}{c|}{min} & \multicolumn{1}{c|}{med} & \multicolumn{1}{c||}{max}  & \multicolumn{1}{c|}{min}  & \multicolumn{1}{c|}{med} & \multicolumn{1}{c||}{max} & \multicolumn{1}{c|}{min}  & \multicolumn{1}{c|}{med} & \multicolumn{1}{c|}{max}  \\ 
\hline
$\Pr\left(\mathcal{M}_{1}\vert y,w=\frac{1}{3}\right)>0.5$ & 16.0 & 25.5 & 10.4 & 79.0 & 74.5 & 47.7 & 90.3 & 89.1 & 63.5 & 96.2 & 93.7 & 74.0 \\
$\Pr\left(\mathcal{M}_{1}\vert y,w=\frac{1}{3}\right)>p$ & 5.0 & 5.0 & 5.0 & 67.6 & 57.9 & 39.3 &  84.7 & 80.1 & 57.8 & 94.0 & 89.0 & 69.2 \\
$\Pr\left(\mathcal{M}_{1}\vert y,w=\frac{1}{3}\right)>0.8$ & 5.7 & 7.4 & 3.4 & 68.5 & 63.9 & 36.2 & 85.2 & 83.5 & 55.1 & 94.2 & 91.6 & 67.5 \\
\hline
$\Pr\left(\mathcal{M}_{1}\vert y,w=\frac{1}{2}\right)>0.5$ & 16.2 & 22.1 & 9.1 & 85.6 & 78.5 & 47.7 & 95.4 & 91.2 & 63.3 & 98.1 & 96.3 & 75.6 \\
$\Pr\left(\mathcal{M}_{1}\vert y,w=\frac{1}{2}\right)>p$ & 5.0 & 5.0 & 5.0 & 75.6 & 64.3 & 38.0 & 92.4 & 84.3 & 56.9 & 97.0 & 93.5 & 70.7 \\
$\Pr\left(\mathcal{M}_{1}\vert y,w=\frac{1}{2}\right)>0.8$ & 3.8 & 5.0 & 2.4 & 73.1 & 64.2 & 30.6 & 91.0 & 84.3 & 49.2 & 96.7 & 93.3 & 66.5  \\
\hline
\end{tabular}
}
\label{tab_simu1_w}
\end{table}

Table \ref{tab_simu_var12_w} is equivalent to Table \ref{tab_simu_var12} of the previous subsection. 
The Bayesian method seemed conservative with the $\mathcal{M}_0$ when the weighting was introduced. 
Thus, we found the same conclusions than the simulations without the weights, but the Bayesian method power was lower with the weight than without.
In light of these results, we advise to consider a higher threshold with the weighting than without. 
The Bayesian method using $p=0.9$ was more powerful than Student test and ensured the good frequentist properties.

\begin{table}[!t]
\caption{Efficiency (in percentage) of methods to detect a difference between the two sample expectations when $(\sigma_{t_i})_{i=1,2}$ changes, given a difference of $\delta$ percent and the weighting $w$. 
$N=1000$ is the simulation number where : $n=10$, $\mu_c=2.94$, $\sigma_{c}=0.6$ et $\sigma_{t_3}=\left(\sigma_{t_k}\right)_{k\neq i}=1.5$. Min corresponds to $\sigma_{t_{i}}=0.6$ and max to $\sigma_{t_{i}}=3$.}
\centering
\scalebox{0.9}{
\begin{tabular}{|c|l||c|c||c|c||c|c||c|c|} 
\hline
\multicolumn{2}{|c|}{Effect} & \multicolumn{2}{c||}{$0\%$} & \multicolumn{2}{c||}{$30\%$} & \multicolumn{2}{c||}{$40\%$} & \multicolumn{2}{c|}{$50\%$} \\ 
\hline
\multicolumn{2}{|c|}{$\sigma_{t_{i}}$}  & \multicolumn{1}{c|}{min}  & \multicolumn{1}{c||}{max} & \multicolumn{1}{c|}{min}  & \multicolumn{1}{c||}{max} & \multicolumn{1}{c|}{min}  & \multicolumn{1}{c||}{max}  & \multicolumn{1}{c|}{min}  & \multicolumn{1}{c|}{max} \\ 
\hline
\hline
$i=1$& $\Pr\left(\mathcal{M}_1\vert y,w=\frac{1}{3}\right)>0.5$ & 2.5 & 5.3 & 30.5 & 45.9 & 45.9 & 70.7 & 63.0 & 87.2 \\
& $\Pr\left(\mathcal{M}_1\vert y,w=\frac{1}{3}\right)>p$ & 5 (p=0.187) & 5 (p=0.515) & 37.8 & 45.7 & 52.5 & 69.5 & 67.6 & 86.7 \\
 &$\Pr\left(\mathcal{M}_1\vert y,w=\frac{1}{3}\right)>0.8$ & 1.1 & 1.7 & 25.7 & 34.2 & 40.0 & 58.7 & 57.1 & 80.7 \\
\hline
$i=2$ &  $\Pr\left(\mathcal{M}_1\vert y,w=\frac{1}{3}\right)>0.5$ & 1.7 &   8.6 & 85.5 & 28.1 & 97.4 & 32.9 & 99.6 & 43.3  \\
&  $\Pr\left(\mathcal{M}_1\vert y,w=\frac{1}{3}\right)>p$ & 5 (p=0.255) & 5 (p=0.76) & 89.6 & 22.5 & 98.5 & 28.4 & 99.9 &  38.7\\
&  $\Pr\left(\mathcal{M}_1\vert y,w=\frac{1}{3}\right)>0.8$ & 0.5 & 3.4 &  77.6 & 19.5 & 95.6 & 26.5 & 99.4 &  35.7 \\
\hline
\hline
$i=1$&  $\Pr\left(\mathcal{M}_1\vert y,w=\frac{1}{2}\right)>0.5$ & 3.2 & 7.8 & 35.6 & 57.3 & 54.3 & 78.6 & 72.6 & 93.0 \\
&  $\Pr\left(\mathcal{M}_1\vert y,w=\frac{1}{2}\right)>p$ & 5 (p=0.367) & 5 (p=0.581) & 39.4 & 52.6 & 58.0 & 75.5 & 75.0 & 91.7\\
 &  $\Pr\left(\mathcal{M}_1\vert y,w=\frac{1}{2}\right)>0.8$ & 0.9 & 1.3 & 26.8 & 35.7 & 46.2 & 63.1 & 64.9 & 83.7 \\
\hline
$i=2$&  $\Pr\left(\mathcal{M}_1\vert y,w=\frac{1}{2}\right)>0.5$ & 3.3 & 8.1 & 92.0 & 33.5 & 98.4 & 41.9 & 99.8 & 52.8\\
&  $\Pr\left(\mathcal{M}_1\vert y,w=\frac{1}{2}\right)>p$ & 5 (p=0.429) & 5 (p=0.604) & 93.1 & 29.9 & 98.9 & 39.5 & 99.8 & 48.4\\
 &  $\Pr\left(\mathcal{M}_1\vert y,w=\frac{1}{2}\right)>0.8$ & 0.9 & 2.4 & 82.6 & 23.0 & 95.9 & 31.9 & 99.7 & 40.7 \\
\hline
\end{tabular}
}
\label{tab_simu_var12_w}
\end{table}

Table \ref{tab_simu_loc_w} shows the trend for the Bayesian approach to detect wrongly the difference when the weighting is taken into account (situations 6 to 9).
Contrary to expectation, the proposed method was influenced more by the prior experiments given the weighting than without.
Indeed, for the situations 1 and 2, the power to detect a difference between the two groups for the third experiment decreased with the weighting (Table \ref{tab_simu_loc_w} versus Table \ref{tab_simu_loc1}). 
In contrast, the situations 3 and 4 showed the power was more important with weighting than without weighting, while the difference between the two groups was mainly observed for the prior experiments.

\begin{table}[!h]
\caption{Efficiency (in percentage) of methods to detect a difference between the two sample expectations, given the weighting $w$, when the location parameter $\left(\mu_{t_i}\right)_{i=1,2,3}$ changes from $\left(\mu_{c_i}\right)_{i=1,2,3}=2.94$ between experiments following a percentage $\delta$.
}
\centering
\begin{tabular}{|l|c|c|c|c|c|c|c|c|c|} 
\hline
\hline
\multicolumn{1}{|c|}{\textbf{Situations}} & \textbf{1} & \textbf{2} & \textbf{3} & \textbf{4} & \textbf{5} & \textbf{6} & \textbf{7} & \textbf{8} & \textbf{9} \\
\hline
\multicolumn{1}{|c|}{$\delta$ of Experiment 1} & $10\%$ & $10\%$ & $30\%$ & $50\%$ & $50\%$ & $50\%$ &$50\%$ & $10\%$ & $10\%$ \\
\multicolumn{1}{|c|}{$\delta$ of Experiment 2} & $10\%$ & $10\%$ & $30\%$ & $30\%$ & $10\%$ & $10\%$ &$30\%$ & $10\%$ & $20\%$ \\
\multicolumn{1}{|c|}{$\delta$ of Experiment 3}  & $40\%$ & $50\%$ & $10\%$ & $10\%$ & $10\%$ & $0\%$ & $0\%$ & $0\%$ & $0\%$ \\
\hline
  $\Pr\left(\mathcal{M}_{1}\vert y,w=\frac{1}{3}\right)>0.5$ & 47.6 &  51.2 & 67.1 & 54.0&  21.7 & 14.5 & 48.7 & 30.9 & 45.0\\
  $\Pr\left(\mathcal{M}_{1}\vert y,w=\frac{1}{3}\right)>0.8$ & 30.4 & 32.8 & 53.7 & 44.7 & 14.2  & 9.9 & 37.5 & 13.9 & 28.8\\
  $\Pr\left(\mathcal{M}_{1}\vert y,w=\frac{1}{3}\right)>0.9$ & 22.7 & 25.8 & 45.2 & 39.7& 11.8 &  8.1 & 31.0 &  9.6 & 20.5 \\
\hline
 $\Pr\left(\mathcal{M}_{1}\vert y,w=\frac{1}{2}\right)>0.5$ & 56.0 & 60.8 & 67.4 & 53.4 & 22.7 & 14.0&  43.7 & 29.6 & 42.3 \\
 $\Pr\left(\mathcal{M}_{1}\vert y,w=\frac{1}{2}\right)>0.8$ & 32.3 & 36.8 & 45.8 & 40.5&  12.9 & 7.1 & 28.9 & 10.5 & 20.1\\
  $\Pr\left(\mathcal{M}_{1}\vert y,w=\frac{1}{2}\right)>0.9$ &  22.3 & 27.4 & 36.3 & 32.7  & 9.3  & 4.8 & 22.0 & 5.3 & 13.0 \\
 \hline
\end{tabular}
\label{tab_simu_loc_w}
\end{table} 

The different simulation studies showed the importance of the choice of experiments to calibrate the prior. The weighting did not seem to clarify or reduce the prior information in our case, and this choice is more delicate with the weighting. 
The choice of these historical experiments is important and they have to be similar (protocol, compound, ...) to the experiment of interest.

\section{Real data analysis}

In this section, the Bayesian approach was applied to experimental results obtained in a model of inflammation, for which historical behavioral data are available.
The objective of the study is to evaluate the effect of test compounds in mice after intra-plantar administration of a pro-inflammatory substance, Complete Freund's Adjuvant (CFA).
CFA rapidly produces local inflammation at the side of administration, and responses measuring spontaneous behaviors (licking, protection of the administrated paw) or behaviors such as paw withdrawal evoked by mechanical of thermal stimuli.
The data from these experiments are measured as reaction time (latency) of mice to withdraw their paw upon stimulation.
The aim is to compare the effects of a compound versus those of a vehicle.

\begin{figure}[h!]
\centering
\includegraphics[scale=0.7]{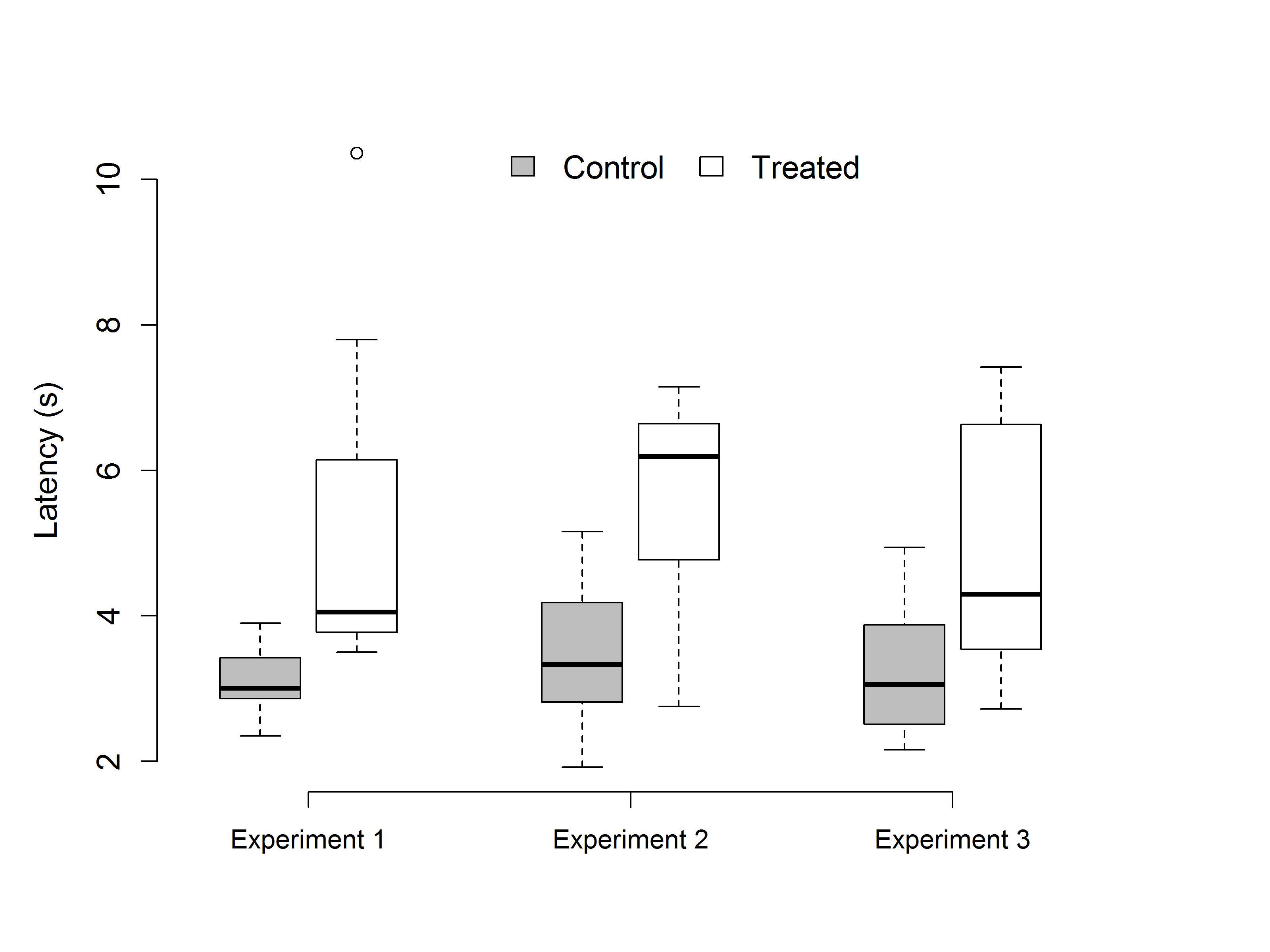}
\caption{Distribution of latency in seconds for each group of the three experiments. Each experiment is composed of two groups: control (mice receiving the vehicle) and treated (mice receiving Ibuprofen). The third experiment (right panel) is the experiment of interest, and the experiments 1 and 2 (left and middle panels) are prior datasets.}
\label{boxplot_application}
\end{figure}

In the present \textit{in vivo} experiments, a radiant heat source was focused on the inflamed hind-paw of every mouse in each group, and the measured variable was a latency to withdraw the paw from the heat source in seconds.
Previous analyses on historical data showed that normality was quite satisfying, but there was often a problem of heterogeneity of variance. 
The selection of historical experiments was made by researchers following our recommendations. 
The three experiments were similar: the control groups received saline solution (vehicle) and Ibuprofen for the treated groups. The sample size was eight for each group. 
Figure \ref{boxplot_application} shows the latency distribution for the three used experiments. The latency was more important in the treated group than in the control group.
In addition to the same conditions, Figure \ref{boxplot_application} shows a similar trend for the three experiments. We are in an ideal case for the application of the three-step Bayesian method.

The p-value associated with the Student test with Satterthwhaite correction was $0.042$ for experiment 3. The null hypothesis was then rejected but the p-value was near to the usual threshold of $5\%$. 
The heterogeneity between variances would have been greater given the lower sample size, the Student test would reject a significant effect of Ibuprofen in comparison with control. 
Our Bayesian method allows to make a statistical decision for the last experiment given the information from the two first experiments.
The Bayesian method advantage is the direct interpretation of the posterior probabilities of models. 
Without assuming a weighting on historical experiments since this was not pertinent based on the simulation results, the posterior probability associated with a similar expectation for the two groups was $\Pr\left(\mathcal{M}_0\vert y\right)=0.006$, and the posterior probability associated with the effect of Ibuprofen compared to the vehicle was $\Pr\left(\mathcal{M}_1\vert y\right)=0.994\%$.
The posterior probability of $\mathcal{M}_1$ reflected the belief in difference between the two groups given the historical data. 
Thus, whatever the threshold $p$ previously discussed to validate the hypothesis $\mathcal{H}_1$ associated with $\Pr\left(\mathcal{M}_1\vert y\right)$, we easily concluded that the latency was different between the two groups.
Both approaches concluded identically for this example: there was an effect of Ibuprofen on the paw withdrawal latency to heat in the CFA protocol. 
The treated group latency was statistically greater than the control group latency. 
Moreover, the Bayesian approach seemed to be a safer method to discriminate between the two hypotheses.


\section{Discussion}

We propose a fully Bayesian technique to solve the Behrens-Fisher problem as an alternative to the usual frequentist method which is not appropriate and could lack of power with small sample size. The corresponding hypothesis testing is considered as a model choice question in the Bayesian paradigm. The proposed method takes into account information from historical data to calibrate an objective informative prior on parameter distribution and also on the model space in two steps. A third and later step is performed to provide the posterior probability of each model. This quantity reflects the confidence in the model (or hypothesis) conditionally on the data. Its interpretation is easy. Despite the small sample sizes, the three step Bayesian method has good frequentist properties as seen in the results of simulation, thanks to the objective informative prior built with the historical data. 

Usually, the chosen model is the one with the greatest posterior probability. Considering two models, the natural threshold of $0.5$ can be used to detect an expectation difference between the two groups ($\mathcal{M}_1$) or not ($\mathcal{M}_0$). The arbitrary choice of this threshold is a difficult task. The threshold of $0.8$ corresponding to our expectation on this type of study is interesting according to our simulation studies because of its good frequentist properties. Optimally, the threshold could be calibrated for each specific protocol in order to control the type I error at $5\%$. 

In most cases, the proposed methodology is more powerful than the others. However, as any Bayesian method, it depends on a pertinent choice of historical data. In the preclinical research context, where the \textit{in vivo} pharmacology experiments are routinely performed using the same protocol, several experiments, which meet the Pocock's criteria \citep{pocock_combination_1976}, could be used to apply the three steps methodology proposed. 
The method naturally weights the two historical datasets and an overweight is not needed. 
Indeed, the historical data used in the second step carries more information (the informative prior probability in the model space is built on the second experiment only). The order of the historical data used is important. We recommend, considering the second step, the most recent historical experiment to decrease potential biases. 
If only one historical experiment is available, the proposed Bayesian method is also suitable in two steps and the results are the posterior probabilities of the second step. However this probability would not be as accurate as it would be after three steps since the informative prior on model space would be missing.

As mentioned in the introduction, the hierarchical models \citep{neuenschwander_summarizing_2010} for the prior building were not studied in this work. However, it could be interesting to expand our work with the comparison between this kind of MAP prior and also with power prior approach \citep{ibrahim_power_2015}. In this paper, we focused on the comparison of two groups. But, in most cases, the objective of this kind of preclinical research experiments is to compare more than one group with a control one, at least three groups are considered: control, treated and reference. 
It would be interesting to develop the proposed method to meet the objective in the case of several groups \citep{scott_exploration_2005}.
With three groups for example, we have to discriminate between four models such as: 
\begin{itemize}
\item Model $\mathcal{M}_{0}$ : no difference between the means of the three groups or no effect of reference compound and tested compound;
\item Model $\mathcal{M}_{1}$ : just reference compound has an effect;
\item Model $\mathcal{M}_{2}$ : just tested compound has an effect;
\item Model $\mathcal{M}_{3}$ : both compounds have an effect versus placebo.
\end{itemize}

\section*{Acknowledgements}

We would like to thank Nicolas Bonnet from Sanofi CEP statistical department in Montpellier for his support and review. Thanks also to Noelle Boussac from Sanofi research statistical department for her support on this work, and thanks to the researchers from Aging Pain Sanofi department in Montpellier for whom these analyses were performed; thanks in particular to Isabel A. Lefevre for her editorial assistance and Veronique Menet, Eric Dufour for providing the experimental data and their availability for exchanges.

\bibliographystyle{apalike}
\bibliography{ArticleBayesianMethod}

\end{document}